\documentstyle[12pt]{article}
\begin{document}
\def\theequation{\thesection.\arabic{equation}}
\newenvironment{proof}{\noindent {\em Proof. }}{\hfill$\Box$.\\}
\newtheorem{theo}{Theorem}[section]
\newtheorem{prop}[theo]{Proposition}
\newtheorem{lemm}[theo]{Lemma}
\newtheorem{coro}[theo]{Corollary}
\newtheorem{defi}[theo]{Definition}
\newtheorem{rema}[theo]{Remark}
\newcommand{\sectio}[1]{\section{#1}\setcounter{equation}{0}}
\newcommand{\e}{\label}
\newcommand{\r}[1]{(\ref{#1})}
\newcommand{\re}{\ref}
\newcommand{\k}{\ldots}
\newcommand{\p}{\partial}
\newcommand{\di}{\p\!\!\!/}
\newcommand{\Bo}{\mbox{\raisebox{-0.3ex}{\Large\mbox{$\Box$}}}}
\newcommand{\api}{{\cal A}}
\newcommand{\bpi}{{\cal B}}
\newcommand{\cpi}{{\cal C}}
\newcommand{\fpi}{{\cal F}}
\newcommand{\lpi}{{\cal L}}
\newcommand{\npi}{{\cal N}}
\newcommand{\ppi}{{\cal P}}
\newcommand{\spi}{{\cal S}}
\newcommand{\sw}{{\varsigma}}
\newcommand{\gpi}{{\cal G}}
\newcommand{\vhi}{\varphi}
\newcommand{\nb}{\mbox{{\bf N}}{}}
\newcommand{\rb}{\mbox{{\bf R}}{}}
\newcommand{\cb}{\mbox{{\bf C}}{}}
\newcommand{\jb}{\mbox{{\bf 1}}{}}
\newcommand{\zb}{\mbox{{\bf Z}}{}}
\newcommand{\xb}{\mbox{{\bf X}}{}}
\newcommand{\cbg}{\cb_{*}}
\newcommand{\rbg}{\rb_{*}}
\newcommand{\al}{\alpha}
\newcommand{\ga}{\gamma}
\newcommand{\del}{\delta}
\newcommand{\eps}{\varepsilon}
\newcommand{\es}{s}
\newcommand{\lam}{\lambda}
\newcommand{\si}{\sigma}
\newcommand{\ta}{\tau}
\newcommand{\de}{\Delta}
\newcommand{\la}{\Lambda}
\newcommand{\wb}{\bar w}
\newcommand{\wc}{{}^cw}
\newcommand{\g}{\Gamma}
\newcommand{\gd}{\dot\Gamma}
\newcommand{\qpg}{quantum Poincar\'e  group}
\newcommand{\qig}{quantum inhomogeneous group}
\newcommand{\be}{\begin{equation}}
\newcommand{\ee}{\end{equation}}
\newcommand{\bt}{\begin{theo}}
\newcommand{\et}{\end{theo}}
\newcommand{\bp}{\begin{prop}}
\newcommand{\ep}{\end{prop}}
\newcommand{\bl}{\begin{lemm}}
\newcommand{\el}{\end{lemm}}
\newcommand{\bc}{\begin{coro}}
\newcommand{\ec}{\end{coro}}
\newcommand{\bde}{\begin{defi}}
\newcommand{\ede}{\end{defi}}
\newcommand{\br}{\begin{rema}}
\newcommand{\er}{\end{rema}}
\newcommand{\bd}{\begin{proof}}
\newcommand{\ed}{\end{proof}}
\newcommand{\ba}{\begin{array}}
\newcommand{\ea}{\end{array}}
\newcommand{\btr}{\begin{trivlist}}
\newcommand{\etr}{\end{trivlist}}
\newcommand{\lra}{\longrightarrow}
\newcommand{\ot}{\otimes}
\newcommand{\op}{\oplus}
\newcommand{\ti}{\times}
\newcommand{\sd}{\rhd\mbox{\hspace{-2ex}}<}              
\newcommand{\spa}{\mbox{{\rm{span}}}}                 
\newcommand{\dm}{\mbox{{\rm{dim}}}}                 
\newcommand{\po}{\mbox{{\rm{Poly}}}}
\newcommand{\id}{\mbox{{\rm{id}}}}
\newcommand{\mor}{\mbox{{\rm{Mor}}}}
\newcommand{\te}{\tilde}
\newcommand{\tb}{\te{\cal B}}
\newcommand{\iffi}{\Leftrightarrow}
\newcommand{\mod}{\mbox{{\rm{mod}}}}
\newcommand{\im}{\mbox{{\rm{im}}}}
\newcommand{\Sp}{\mbox{{\rm{Sp}}}}
\newcommand{\tr}{\mbox{{\rm{tr}}}}
\newcommand{\Irr}{\mbox{{\rm{Irr}}}}
\newcommand{\Rep}{\mbox{{\rm{Rep}}}}
\newcommand{\diag}{\mbox{{\rm{diag}}}}
\newcommand{\ov}{\overline}
\newcommand{\gw}{${}^*$}
\newcommand{\ite}{\item[]}
\newcommand{\hs}{\hspace{0.2ex}}
\newcommand{\pri}{{}^{\prime}}
\newcommand{\itemite}{\item[]\hspace{3ex}}
\newcommand{\Ren}{\mbox{{\rm{Re}}}}
\newcommand{\Imn}{\mbox{{\rm{Im}}}}

\title{The Dirac operator and gamma matrices\\ 
for quantum Minkowski spaces}
\author{P. Podle\'s\thanks{This research was supported in part
by NSF grant DMS--9508597 and in part by Polish KBN grant No. 2 P301 020 07}\\
Department of Mathematical Methods in Physics\\
Faculty of Physics, University of Warsaw\\
Ho\.za 74, 00--682 Warszawa, Poland}
\date{March 6, 1997}
\maketitle
\begin{abstract}
Gamma matrices for quantum Minkowski spaces are found. The invariance of the 
corresponding Dirac operator is proven. We introduce momenta for
spin $1/2$ particles and get (in certain cases) formal solutions
of the Dirac equation.
\end{abstract}
\setcounter{section}{-1}
\sectio{Introduction}

It is widely recognized that geometry of space--time should
drastically change at very small distances, comparable with
Planck's length. On the other hand there is no satisfactory
physical theory which would describe such a change. A lot of
effort was devoted to simple physical models describing possible
changes of geometry which can occur. One of possibilities is
provided by the theory of quantum groups: the related examples 
of quantum space--times have a (quantum) group of symmetries
which can be as big as the classical Poincar\'e group. If we
want to have a quantum space--time which has exactly the same
properties as the classical Minkowski space endowed with the
action of (spinorial) Poincar\'e group, we get the
classification of quantum Minkowski spaces and quantum
Poincar\'e groups given in \cite{POI}. The related differential
structure on quantum Minkowski spaces was determined in
\cite{KG}. Therefore we are able to write Klein--Gordon equation
for spin $0$ particles and solve it (at least formally) in many
cases. The same holds for Dirac equation for spin $1/2$
particles provided we are able to find the gamma matrices (and
certain other objects). This remaining task is undertaken in the
present paper. In Section~1 we recall the theory of quantum
Minkowski spaces and \qpg s. In Section~2 we prove that the
requirement of invariance of the Dirac operator determines all
the gamma matrices up to two constants $a,b\in\cb$. The square of the Dirac
operator is equal to the Laplacian (as in the classical case) if
and only if $ab=1$. The normalization $a=b=1$ is chosen.
We study certain expressions like the deformed Lagrangian.
In Section~3 we get (in certain cases) formal solutions of the
Dirac equation and introduce
the momenta for spin $1/2$ particles.

The gamma matrices for
the example of the standard $q$-Lorentz group and $q$-Minkowski
space (\cite{PW},\cite{CSSW}) were considered in
\cite{CSSW} (cf. also \cite{So}, \cite{Sch},\cite{Az}).
This case however doesn't fall
into our scheme (the corresponding $q$-Poincar\'e group contains
dilatations) and involves essentially only one $2$-dimensional
$R$-matrix (while in our considerations 
we have two independent $2$-dimensional $R$-matrices: $L$ and $X$).

The gamma matrices (given by \r{3.6} and \r{2.12}) and relations
among them (the condition 2. of Theorem \re{t2.2} and \r{2.4}) 
were announced in \cite{PP}.
Recently, when the present paper was essentially completed,
some gamma matrices (satisfying a condition like the condition 2. 
of Theorem \re{t2.2}) appeared also in \cite{AA}. That paper
contains explicit formulae for the gamma matrices and metric
tensor in some cases. Also the transformation properties are
discussed there.

We sum over repeated indices (Einstein's convention).  
If
$V,W$ are vector spaces then $\tau: V \otimes W \to W \otimes V$ is given by
$\tau(x \otimes y) = y \otimes x$, $x \in V$, $y \in W$.  We denote the unit
$N\times N$ matrix by $\jb_N$ or $\jb$.
If ${\cal A}$ is an algebra, $v \in M_N({\cal A})$, $w \in
M_K({\cal A})$, then the tensor product $v \otimes w \in M_{NK}({\cal A})$ is
defined by
\[
(v \otimes w)^{ij}{}_{kl} = v^{i}{}_{k}w^{j}{}_{l},\ i,k = 1,\dots,N,\ %
j,l = 1,\dots,K.
\]
We set $\dm\ v = N$.  If ${\cal A}$ is a ${}^*$-algebra then
the conjugate of $v$ is defined as ${\bar v} \in M_N({\cal A})$ where ${\bar
v}^{i}{}_{j} = (v^{i}{}_{j})^*$.  We also set $v^* = {\bar v}^T$ ($v^T$ 
denotes the
transpose of $v$, i.e. $(v^T)_{i}{}^{j} = v^{j}{}_{i}$). We
write $a\sim b$ if $a,b$ are proportional, i.e. $a=kb$ for $k\in\rbg$.

Throughout the paper quantum groups $H$ are abstract objects described by the
corresponding Hopf \gw-algebras $\mbox{Poly}(H) = ({\cal A},\de)$.  We 
denote
by $\de,\eps,S$ the comultiplication, the counit and the coinverse of
$\mbox{Poly}(H)$.  
In particular, $S$ is invertible ($S^{-1}=*\circ S\circ*$).
We say that $v$ is a representation of $H$ ($v\in\mbox{ Rep } H$)
if $v \in M_N({\cal A})$, $N \in {\nb}$, and
\[
\de v^{i}{}_{j} = v^{i}{}_{k} \otimes v^{k}{}_{j},\ \eps(v^{i}{}_{j}) = 
\del^{i}{}_{j},\ i,j = 1,\dots,N,
\]
in which case
 $S(v^{i}{}_{j})=(v^{-1})^{i}{}_{j}$.
The conjugate of a representation and tensor products of representations are
also representations.  
If $v,w \in \mbox{ Rep } H$, $\dm\ v=N$, $\dm\ w=M$, then we say
that $A \in M_{M \ti N}({\cb})$ intertwines $v$ with $w$ 
if $Av = wA$. We say  that $v,w$ are
equivalent ($v\simeq w$) if such $A$ can be chosen as invertible.
For $\rho\in {\cal A}'$ (the dual
vector space of ${\cal A}$), 
$a\in {\cal A}$, we set
$\rho * a=( \id \otimes\rho)\de a$, $a*\rho=(\rho\otimes \id )\de a$.

\sectio{Quantum Lorentz and Poincar\'e groups}

In this section we recall the definitions and properties of
quantum Lorentz and Poincar\'e groups as well as quantum
Minkowski spaces. These objects are the natural generalizations
of the standard objects known from the relativistic physics. 

Quantum Lorentz groups are defined as quantum groups with the
same properties as the classical Lorentz group $SL(2,\cb)$ \cite{WZ}.
The classification of quantum Lorentz groups is given as follows
\cite{WZ}. The set $\api$ of polynomials on a quantum Lorentz
group is the universal unital ${}^*$-algebra generated by
$w_{AB}$, $A,B=1,2$, satisfying
\be w^{A}{}_{B}w^{C}{}_{D}E^{BD}=E^{AC}, \e{1.1}\ee
\be E'_{AC}w^{A}{}_{B}w^{C}{}_{D}=E'_{BD}, \e{1.2}\ee
\be X^{AB}{}_{CD}w^{C}{}_{K}w^{D}{}_{L}{}^*=
w^{A}{}_{C}{}^*w^{B}{}_{D}X^{CD}{}_{KL}, \e{1.3}\ee
$A,B,C,D,K,L=1,2$, where the matrices $E\in M_{4\ti 1}(\cb)$, $E'\in
M_{1\ti 4}(\cb)$, $X\in M_{4\ti 4}(\cb)$ are given 
(up to a nonzero factor) in \cite{WZ}.
The set $\api$ becomes a Hopf ${}^*$-algebra if we define the
comultiplication and the counit in such a way that $w$ becomes a
representation, i.e.
\be \de w^{A}{}_{B}=w^{A}{}_{C}\ot w^{C}{}_{B},\quad
\eps(w^{A}{}_{B})=\del^{A}{}_{B},\quad A,B=1,2. \e{1.4}\ee
The equations \r{1.1}-\r{1.3} can be also written as
\be (w\ot w)E=E,\ \  E'(w\ot w)=E',\ \  X(w\ot\wb)=(\wb\ot w)X. \e{1.4'}\ee
In particular, setting $E^{11}=E^{22}=0$, $E^{12}=1$, $E^{21}=-1$,
$E'_{AB}=-E^{AB}$, $X^{AB}{}_{CD}=\del^{A}{}_{D}\del^{B}{}_{C}$, we get the
classical Lorentz group $SL(2,\cb)$. Then
$w^{A}{}_{B}$ are the matrix elements of the fundamental
representation of $SL(2,\cb)$, i.e. 
$w^{A}{}_{B}(h)=h^{A}{}_{B}\in\cb$, $h\in SL(2,\cb)$. Moreover,
$f^*(g)=\ov{f(g)}\in\cb$ for $f\in\api$, $g\in SL(2,\cb)$. 

For any quantum Lorentz group we define its representation $\la$
as 
\be \la=V^{-1}(w\ot\wb)V \e{1.4a}\ee
where
\be V^{AB}{}_{i}=(\si_i)_{AB},\quad 
(V^{-1})^{i}{}_{AB}=\frac12\ov{(\si_i)_{AB}}
=\frac12(\si_i)_{BA},\e{1.5}\ee
\[ \si_0=\left(\ba{rr} 1&0\\0&1\ea\right),\ \ 
\si_1=\left(\ba{rr} 0&1\\1&0 \ea\right), \ \ 
\si_2=\left(\ba{rr} 0&-i\\i&0\ea\right), \ \ 
\si_3=\left(\ba{rr} 1&0\\0&-1\ea\right)  \]
are the Pauli matrices. Then $\bar\la=\la$ is irreducible. 

In the next step we introduce the \qpg s, i.e. the quantum
groups with the properties of the (spinorial) Poincar\'e group.
Their definition and (almost complete) classification are given
in \cite{POI}. It turns out that each \qpg\ is related to
one of quantum Lorentz groups described by
$E$, $E'$ and $X=\ta Q'$ as follows:
\btr
\item[1)] $E=e_1\ot e_2-e_2\ot e_1$, $E'=-e^1\ot e^2+e^2\ot e^1$,
\[ Q'=\left( \ba{cccc} t^{-1} & 0 & 0 & 0 \\ 0 & t & 0 & 0\\ 0 & 0 & t & 0\\
             0 & 0 & 0 & t^{-1}\ea\right), \quad 0<t\leq1,\quad \mbox{or} \]
\item[2)] \[ E,E' \mbox{ as above },\quad            
 Q'=\left( \ba{cccc} 1 & 0 & 0 & c^2 \\ 0 & 1 & 0 & 0\\ 0 & 0 & 1 & 0\\
             0 & 0 & 0 & 1\ea\right),\quad \mbox{or} \]
\item[3)] $E=e_1\ot e_2-e_2\ot e_1+ce_1\ot e_1$, $E'=-e^1\ot e^2+e^2\ot
e^1+ce^2\ot e^2$,            
 \[ Q'=\left( \ba{cccc} 1 & 0 & 0 & rc^2 \\ 0 & 1 & 0 & 0\\ 0 & 0 & 1 & 0\\
             0 & 0 & 0 & 1\ea\right),\quad r\geq0,\quad \mbox{ or } \]
\item[4)] \[ E,E' \mbox{ as above },\quad
 Q'=\left( \ba{cccc} 1 & c & c & 0 \\ 0 & 1 & 0 & -c\\ 0 & 0 & 1 & -c\\
             0 & 0 & 0 & 1\ea\right),\quad  \mbox{ or } \]
\item[5)] $E=e_1\ot e_2+e_2\ot e_1$, $E'=e^1\ot e^2+e^2\ot e^1$,  
 \[ Q'=i\left( \ba{cccc} t^{-1} & 0 & 0 & 0 \\ 0 & -t & 0 & 0\\ 
 0 & 0 & -t & 0\\
             0 & 0 & 0 & t^{-1}\ea\right),\quad 0<t\leq1,\quad 
             \mbox{ or } \]
\item[6)] \[ E,E' \mbox{ as above },                 
 Q'=i\left( \ba{cccc} 1 & 0 & 0 & c^2 \\ 0 & -1 & 0 & 0\\ 0 & 0 & -1 & 0\\
             0 & 0 & 0 & 1\ea\right),\quad \mbox{ or } \]
\item[7)] \[ E,E' \mbox{ as above },\quad
 Q'=i\left( \ba{cccc} r & 0 & 0 & \sw \\ 0 & -r & \sw & 0\\ 0 & \sw & -r & 0\\
             \sw & 0 & 0 & r\ea\right), \]             
\[ r=(t+t^{-1})/2,\quad \sw=(t-t^{-1})/2,\quad 0<t<1, \]
\etr
$e_1,e_2$ form the standard basis of $\cb^2$, $e^1,e^2$
the corresponding dual basis of $(\cb^2)^*$, $c\in\rbg$ (cf.
Remark below). One has $q=1$ in
the cases 1)--4) and $q=-1$ for 5)--7). 
We set $q^{1/2}=1$ for $q=1$ and $q^{1/2}=i$ for $q=-1$.

The set $\bpi$ of
polynomials on a \qpg\ is the universal unital \gw-algebra
generated by
$\api$ and $y^i$, $i=0,1,2,3$, satisfying $I_{\bpi}=I_{\api}$,
\be(R-\jb)^{ij}{}_{kl}(y^ky^l-Z^{kl}{}_{m}y^m+
T^{kl}-\la^{k}{}_{m}\la^{l}{}_{n}T^{mn})=0,\e{1.6}\ee
\be y^iw^{A}{}_{B}=G^{iA}{}_{Cj}w^{C}{}_{B}y^j+(H_V)^{iA}{}_{C}w^{C}{}_{B}-
\la^{i}{}_{j}w^{A}{}_{C}(H_V)^{jC}{}_{B},
\e{1.7}\ee
\be (y^{i})^{*}=y^i,\e{1.8}\ee
where 
$R=(V^{-1}\ot V^{-1})(\jb_2\ot X\ot\jb_2)(L\ot\te L)(\jb_2\ot X^{-1}\ot\jb_2)
(V\ot V)$,
$G=(V^{-1}\ot\jb_2)(\jb_2\ot X)(L\ot\jb_2)(\jb_2\ot V)$,
$Z=(\jb_4\ot V^{-1})[H_V\ot\jb_2+(G\ot\jb_2)(\jb_2\ot\te H_V)]V$,
$\te H_V=-\ta\ov{G^{-1}H_V}$, 
$L=\es q^{1/2}(\jb_4+qEE')$, $\te L=q\tau L\tau$,
$s=\pm1$, $(H_V)^{iC}{}_{D}=(V^{-1})^{i}{}_{AB}H^{ABC}{}_{D}$, 
$T^{ij}=(V^{-1})^{i}{}_{AB}(V^{-1})^{j}{}_{CD}T^{ABCD}$,
the possible $H^{ABC}{}_{D},T^{ABCD}\in\cb$ (for given quantum
Lorentz group and $s$), $A,B,C,D=1,2$, are provided in
\cite{POI} (for $c,k=1$).
Then $\bpi$ becomes a Hopf \gw-algebra if we define the
comultiplication and counit in such a way that
\be \de y^i=\la^{i}{}_{j}\ot y^j+y^i\ot I,\quad \eps(y^i)=0,
\quad i=0,1,2,3,\e{1.9}
\ee
and \r{1.4} holds.

{\bf Remark.} We set $c=k=1$ for the objects considered in
\cite{POI}. The objects with $c,k\in\rbg$ 
($c=1$ in cases 1)5)7)) are isomorphic to
those with $c=k=1$ in the following way: $N=\diag(c,1)$,
$\npi=V^{-1}(N\ot N)V$, $w(c)=NwN^{-1}$, $X(c)=(N\ot
N)X(N^{-1}\ot N^{-1})$, $E(c)\sim(N\ot N)E$, $E'(c)\sim
E'(N^{-1}\ot N^{-1})$, $\la(c)=\npi\la\npi^{-1}$, $y^i(c,k)=k\npi^i{}_jy^j$,
$H(c,k)=k(N\ot N\ot N)HN^{-1}$, $T(c,k)=k^2(N\ot N\ot
N\ot N)T$ etc. One should mention that formally taking
$k\lra 0$,
next $c\lra 0$ (in cases 2)3)4))
or $t\lra 1$ (in case 1)), we can deform
all the objects related to the cases 1)-4) with $s=1$ (and any
allowed $H,T$) to the classical case
[ 1), $s=t=1$, $H=T=0$]. Therefore $k,c$ and $1-t$
can serve as small deformation parameters.\\

In particular, $H=T=0$ is always allowed. If in addition we
consider the classical Lorentz group and $s=1$, we get the
classical (spinorial) Poincar\'e group. It is defined as
\[ P=SL(2,\cb)\sd\rb^4=\{(g,a):g\in SL(2,\cb),a\in\rb^4\} \]
with multiplication $(g,a)\cdot(g',a')=(gg',a+\lam_g(a'))$ where the double
covering $SL(2,\cb)\ni g\lra \lam_g\in SO_0(1,3)$ is given by 
$\lam_g(x)^i\si_i=g(x^j\si_j)g^*$, $g\in SL(2,\cb)$, $x\in\rb^4$
(we treat $\lam_g$ as mapping from $\rb^4$ into $\rb^4$). The group
$P$ is the double covering of the (connected component of)
vectorial Poincar\'e group
\[ \te P=SO_0(1,3)\sd\rb^4=\{(M,a): M\in SO_0(1,3), a\in\rb^4\} \]
with multiplication $(M,a)\cdot(M',a')=(MM',a+Ma')$.
This covering  $\pi:P\lra\te P$ is 
defined by $\pi(g,a)=(\lam_g,a)$. We should mention that $P$ is
more important in quantum field theory than $\te P$. In these notations
$f(g,a)=f(g)$ for $f\in\api$ (in particular,
$w^{A}{}_{B}(g,a)=w^{A}{}_{B}(g)=g^{A}{}_{B}\in\cb$) and $y^i(g,a)=a^i=i$-th
coordinate of $a\in\rb^4$. Then the relations \r{1.1}--\r{1.3}
and \r{1.6}--\r{1.7} express the commutativity of our algebra
($R=\tau$, $G=\tau$, $Z=0$ in this case). The case 1), $s=1$, $0<t<1$,
$H=T=0$, corresponds to the \qpg\ of \cite{Ch}. 

Finally, we pass to quantum Minkowski spaces, i.e. the quantum
analogues of the Minkowski space. Their definition and
properties are also given in \cite{POI}. According to \cite{POI},
each \qpg\ admits exactly one quantum Minkowski space. 
The set $\cpi$ of polynomials on a quantum Minkowski space is
introduced as the universal unital ${}^*$-algebra generated by $x^i$, 
$i=0,1,2,3$, satisfying the relations
\be (R-\jb_4)^{ij}{}_{kl}(x^k x^l-Z^{kl}{}_{m}x^m+T^{kl})=0,\e{1.10}\ee
\be (x^{i})^{*}=x^i, \e{1.11}\ee
$i,j=0,1,2,3$. The action of the \qpg\ on the quantum Minkowski
space is given by the unital \gw-homomorphism
$\Psi :\cpi\lra \bpi\ot\cpi$ satisfying
$(id\ot\Psi)\Psi=(\Delta\ot\id)\Psi$, $(\eps\ot\id)\Psi=\id$ and
\be \Psi(x^i)=\la^{i}{}_{j}\ot x^j+y^i\ot I,\quad i=0,1,2,3.\e{1.12}\ee

In particular, for the classical Poincar\'e group we get the classical
Minkowski space $M=\rb^4$. Then $x^i$ are the coordinates on
$M$: $x^i(v)=v^i=i$-th coordinate of $v\in M$. Moreover, $\Psi$
corresponds to the action $\si:P\ti M\lra M$ of $P$ on $M$ given
by $\si(p,v)=(g,a)v=\lam_g(v)+a$, $p=(g,a)\in P$, $v\in M$.
One has
\be \Psi(f)=f\circ\si \e{1.12a}\ee
for $f\in\cpi$. 

Let us recall \cite{KG} that $4$-dimensional covariant
differential calculus on a quantum Minkowski space exists if and
only if 
\[ \te F\equiv[(R-\jb)\ot\jb_4]\{(\jb_4\ot Z)Z-(Z\ot\jb_4)Z+T\ot\jb_4-
(\jb_4\ot R)(R\ot\jb_4)(\jb_4\ot T)\}=0. \]
This requirement singles out some \qpg s
which are described after the proof of Theorem 1.1 of \cite{KG}
(in particular, the trivial choice $H=T=0$ is always allowed). 
{}From now on we limit ourselves to \qpg s and quantum Minkowski
spaces with $\te F=0$. Then for given quantum Minkowski space
the $4$-dimensional covariant differential calculus exists and is unique.
It is described by partial derivatives $\p_i:\cpi\lra\cpi$, 
$i=0,1,2,3$, which are determined by the following properties:
\be \p_i(I)=0,\e{1.13}\ee
\be \p_i(x^kf) = \del^{k}{}_{i}f + (R^{kl}{}_{in} x^n +
Z^{kl}{}_{i})(\p_lf),\quad f\in\cpi. \e{1.14}\ee
The $\p_i$ satisfy the following covariance property
\be (\id \otimes \p_j)(\Psi f) = (\Lambda^{i}{}_{j} \otimes I)
[\Psi(\p_if)],\quad f\in\cpi. \e{1.15}\ee

According to Proposition 3.1.2 of \cite{KG}, its proof and \cite{POI},
$\p_i$ can be also obtained as follows. We set
$\te G=(V^{-1}\ot\jb_2)(\jb_2\ot\te L)(X^{-1}\ot\jb_2)(\jb_2\ot V)$
and define a unital homomorphism
$f=(f^{i}{}_{j})_{i,j=0}^3:\api\lra M_4(\cb)$ by
$f^{i}{}_{j}(w^{C}{}_{D})=G^{iC}{}_{Dj}$, $f^{i}{}_{j}(w^{C}{}_D{}^*)
=\te G^{iC}{}_{Dj}$.
Then there exists a unital homomorphism $\xb:\bpi\lra M_5(\cb)$
given by

\be \xb(a) = \left( \begin{array}{cc}
(f^{l}{}_{j}(S(a)))_{j,l=0}^3 & 0 \\
0 & \eps(a)
\end{array} \right),\ \ a \in {\cal A},\e{1.15a}\ee
\be \xb(y^i) = \left( \begin{array}{cc}
(Z^{il}{}_{j})_{j,l=0}^3 & (\del^{i}{}_{j})_{j=0}^3 \\
0 & 0
\end{array} \right),\ \ i = 0,1,2,3. \e{1.15b}
\ee
In practical computations one can use the formula \cite{POI} $f^l{}_j(S(a))=
\overline{f^l{}_j(a^*)}$.
It turns out that
\be
\xb = \left( \begin{array}{cc}
(X_{j}{}^{l})_{j,l=0}^3 & (Y_j)_{j=0}^3 \\
0 & \eps
\end{array} \right) \e{1.15b'}
\ee
with $X_j{}^l,Y_j \in {\cal B}'$.  One has
\be
\p_j = (Y_j \otimes \id)\Psi,\ \ 
\rho_{j}{}^{k} = (X_{j}{}^{k} \otimes \id)\Psi \e{1.15c}\ee
where $\rho_{j}{}^{k}:\cpi\lra\cpi$ also appear in the proof of
Theorem~1.1 of \cite{KG}. 

The metric tensor is defined \cite{KG} as an invertible matrix 
$g=(g^{ij})_{i,j=0}^3\in M_{4\ti 4}(\cb)$ such that its matrix
entries $g^{ij}$ satisfy the invariance and self--adjointness conditions:
\be \la^{i}{}_{j}\la^{k}{}_{l}g^{jl}=g^{ik}, \e{1.16}\ee
\[ \ov{g^{ik}}=g^{ki}, \]
$i,k=0,1,2,3$. Such $g$ is unique up to a nonzero real
multiplicative factor. It satisfies $Rg=g$. Here we choose it as
\be g=-2q^{1/2}(V^{-1}\ot V^{-1})(\jb_2\ot X\ot\jb_2)(E\ot\ta E)
\e{1.17}\ee
(cf. \cite{KG} and Remark below). 
Matrix elements of $g^{-1}$ are denoted by $g_{ij}$. 
The Laplacian is given
\cite{KG} as 
\be \Bo=g^{ij}\p_j\p_i. \e{1.18}\ee
It is invariant and commutes with partial derivatives:
\be (\id\ot\Bo)[\Psi(f)]=\Psi(\Bo f),\quad f\in\cpi, \e{1.19}\ee
\be \Bo\p_i=\p_i\Bo,\quad i=0,1,2,3. \e{1.20}\ee
According to \cite{POI} and \cite{KG}, the ``sizes'' of all our
constructions are the same as for the classical Poincar\'e group,
the classical Minkowski space and the standard differential
calculus on it (we consider only the polynomial functions).

{\bf Remark.} The factor (-2) in \r{1.17} gives the standard
metric tensor $g=\diag(1,-1,-1,-1)$ for the classical Poincar\'e
group and was taken into account in the expression for a
propagator near the end of Section~4 of \cite{KG}. However, it
was not fixed in the considerations after (3.8) of \cite{KG}.\\

In particular cases the metric tensor $g=(g^{ij})_{i,j=0}^3$ equals
\btr
\item[1)5)] $g=q\diag(t,-t^{-1},-t^{-1},-t)$,
\item[2)6)] $g=\frac12q\left(\ba{cccc}
2-c^2&0&0&-c^2\\
0&-2&0&0\\
0&0&-2&0\\
-c^2&0&0&-2-c^2\ea\right)$,
\item[3)] $g=-\frac12\left(\ba{cccc}
(r+1)c^2-2&0&-2ic&(1+r)c^2\\
0&2&0&0\\
2ic&0&2&2ic\\
(1+r)c^2&0&-2ic&(r+1)c^2+2\ea\right)$,
\item[4)] $g=-\frac12\left(\ba{cccc}
3c^2-2&-2c&-2ic&3c^2\\
-2c&2&0&-2c\\
2ic&0&2&2ic\\
3c^2&-2c&-2ic&3c^2+2\ea\right)$,
\item[7)] $g=\diag(-r+\sw,r+\sw,r-\sw,r+\sw)$
\etr
($q=1$ for 1)--4) and $q=-1$ for 5)--7)).

\sectio{Invariance of the Dirac operator}

In this section we prove that the requirement of invariance of
the Dirac operator $\di$ determines all the gamma matrices up to
two constants. Then using the condition $\di^2=\Bo$ we provide
the exact form of the gamma matrices. We also study
certain expressions like the deformed Lagrangian.

We shall consider gamma matrices $\ga^i\in M_{4\ti 4}(\cb)$,
$i=0,1,2,3$. At the moment they are not determined yet. The
Dirac operator has form $\di=\ga^i\ot\p_i$ (cf. \cite{KG}). It acts
on bispinor functions $\phi\in\te\cpi\equiv\cb^4\ot\cpi$ (in a more
advanced approach we should consider square integrable functions
$\phi$). 

In the classical case the Poincar\'e group $P$ acts on $\te\cpi$
as follows:
\[ \phi'(x')=\spi(g)\phi(x) \]
where $\phi\in\te\cpi$, $x\in M$, $\phi(x)\in\cb^4$, $\phi'$ is
$\phi$ transformed by $p=(g,a)\in P$, $x'=p\cdot x$,
$\phi'(x')\in\cb^4$, $\spi$ is a representation of the Lorentz
group $SL(2,\cb)$ acting in the space $\cb^4=\cb^2\op\cb^2$ of
bispinors, $\spi\simeq w\op\wb$ (undotted and dotted spinors).
Writing $\phi=\eps_a\ot\phi^a$ where $\eps_a$, $a=1,2,3,4$, form
the standard basis of $\cb^4$, replacing $x$ by $p^{-1}x$ and
then $p$ by $p^{-1}$, using \r{1.12a} and setting
\be [\te\Psi(\phi)](p,x)=\phi'(x),\e{2.0}\ee
$\Psi(\phi^a)=\phi^{a(1)}\ot\phi^{a(2)}$
(Sweedler's notation, exception of Einstein's convention), one obtains
\be \te\Psi(\eps_a\ot\phi^a)=\gpi_{a}{}^{l}\phi^{a(1)}\ot\eps_l\ot\phi^{a(2)}
\e{2.1}\ee
where $\gpi=(\spi^{-1})^T\simeq w\op\wb$ (cf. remarks after \r{2.2}). 

For the general quantum Minkowski space we define the action
$\te\Psi:\te\cpi\lra\bpi\ot\te\cpi$ of \qpg\ on $\te\cpi$ by \r{2.1} where
$\gpi$ is a representation equivalent to $w\op\wb$. For the future
convenience we choose
\be \gpi={}^cw\op\wb \e{2.2}\ee
where ${}^cw=(w^T)^{-1}$. We notice that \r{1.4'} implies (here
we treat $E$ as $2\ti 2$ matrix and then \cite{WZ} $E'=E^{-1}$)
\be wEw^T=E. \e{2.2'} \ee
Therefore ${}^cw\simeq w\ \ $: 
\be {}^cw=E^{-1}wE. \e{2.2a}\ee

In the classical case one has the invariance of the Dirac operator:
\be [\di\phi]'=\di\phi'. \e{2.3}\ee
Then assuming that $\phi$ satisfies the Dirac equation
$(i\di-m)\phi=0$, we get
\[ 0=[(i\di-m)\phi]'=i(\di\phi)'-m\phi'=(i\di-m)\phi',\]
i.e. $\phi'$ also satisfies the Dirac equation. It means that
the Dirac equation is invariant. On the other hand using
\r{2.0}, the invariance condition \r{2.3} is equivalent to
\be \te\Psi(\di\phi)=(\id\ot\di)[\te\Psi(\phi)],\quad \phi\in\te\cpi.
\e{2.4}\ee

We set \r{2.4} as the condition of invariance of the Dirac
operator for the general quantum Minkowski space. The main
result of the present Section is contained in

\bt\e{t2.1} The following are equivalent:\\
\btr
\item[1.] The Dirac operator is invariant (i.e. \r{2.4} is satisfied)
\item[2.] 
\be \ga^i=\left(\ba{cc} 0&bA_i\\
                a\si_i&0\ea\right) \e{2.11}\ee
\etr
where
\be A_i=q^{-1/2}E^T(\si_i\circ D)E,\e{2.12}\ee
$(\si_i\circ D)_{KL}=(\si_i)_{AB}D^{AB}{}_{KL}$, $D=\tau X^{-1}\tau$,
$a,b\in\cb$
($E$ is regarded here as $2\ti 2$ matrix).\et

\bd Setting $\phi=\eps_a\ot\phi^a$ and using \r{2.1}, one easily
checks that the LHS of \r{2.4} equals 
\be
\gpi_{m}{}^{s}(\ga^j)^{m}{}_{a}(\p_j\phi^a)^{(1)}\ot\eps_s
\ot(\p_j\phi^a)^{(2)}.
\e{2.5}\ee
Similarly, using \r{2.1} and \r{1.15}, we obtain that the RHS of
\r{2.4} is equal to
\be \gpi_{a}{}^{l}\phi^{a(1)}\ot(\ga^i)^{s}{}_{l}\eps_s\ot\p_i[\phi^{a(2)}]
=\gpi_{a}{}^{l}\la^{j}{}_{i}(\p_j\phi^a)^{(1)}(\ga^i)^{s}{}_{l}\ot\eps_s
\ot(\p_j\phi^a)^{(2)}.
\e{2.6}\ee
Choosing $\phi=x^r\eps_t$, one gets $\phi^a=x^r\del^{a}{}_{t}$,
$\p_j\phi^a=\del^{r}{}_{j}\del^{a}{}_{t}I$. In this case \r{2.4} (i.e. the
equality of \r{2.5} and \r{2.6}) yields
\be \gpi_{m}{}^{s}(\ga^r)^{m}{}_{t}=\gpi_{t}{}^{l}\la^{r}{}_{i}(\ga^i)^s{}_l.
\e{2.7}\ee
On the other hand \r{2.7} implies that \r{2.5} equals \r{2.6}
and \r{2.4} follows (for general $\phi$). Therefore \r{2.4} and
\r{2.7} are equivalent. 

Setting 
\be W_{l}{}^{i}{}_{,}{}^{s}=(\ga^i)^{s}{}_{l}, \e{2.7a}\ee
we can write \r{2.7} as
\be W\gpi=(\gpi\ot\la)W.\e{2.8}\ee
Defining
\be N=(\jb_4\ot V)W,\e{2.8a}\ee
using \r{1.4a} and \r{2.2}, \r{2.8} can be translated as
\be N({}^cw\op\wb)=({}^cw\ot w\ot\wb\op\wb\ot w\ot\wb)N.\e{2.9}\ee
According to \r{2.2a} and Proposition 2.1 of \cite{POI} (cf.
also \cite{WZ}), 
\[ {}^cw\ot w\ot\wb\simeq(w^1\ot\wb)\op\wb,\quad \wb\ot
w\ot\wb\simeq (w\ot\ov{w^1})\op w \]
(decompositions into irreducible components, we use the notation
of Proposition 2.1 of \cite{POI}). Thus
\be N=\left(\ba{cc} 0&N_1\\N_2&0\ea\right), \e{2.10}\ee
\r{2.9} means that 
\be N_1\wb=({}^cw\ot w\ot\wb)N_1,  \e{2.11'}\ee
\be N_2{}^cw=(\wb\ot w\ot\wb)N_2      \e{2.12'}\ee
and $N_1,N_2$ are fixed up to multiplicative constants. Using
the definition of $\wc$, \r{1.3} and conjugated \r{1.1}, one can
check that solutions of \r{2.11'},\r{2.12'} are given by
\[ (N_1)^{ABC}{}_{D}=2a\del^{A}{}_B\del^{C}{}_{D}, \]
\[ (N_2)^{ABC}{}_{D}=2bq^{-1/2}(X^{-1})^{BC}{}_{KL}E^{KA}E^{LD}, \]
$a,b\in\cb$, where additional scalar factors $2,q^{-1/2}$ are added
for future convenience. Using \r{2.7a}, \r{2.8a} and \r{1.5}, we
finally get \r{2.11}.\ed

In the standard Dirac theory $\di^2=\Bo$. We set this as the
additional condition for our gamma matrices. Then the Dirac
equation $(i\di-m)\phi=0$ implies (formally) the Klein--Gordon
equation $(\Bo+m^2)\phi=0$. 

\bt\e{t2.2} Assume \r{2.11}. The following are equivalent:
\btr
\item[1.] $\di^2=\Bo$.
\item[2.] $\ga^i\ga^j+R^{ji}{}_{lk}\ga^k\ga^l=2g^{ji}\jb,\ \ \ \ \
i,j=0,1,2,3$. 
\item[3.] $ab=1$.
\etr\et

{\bf Remark.} The condition 2 was considered in \cite{KG}
(cf. \cite{RTF}, \cite{Sch}). For the classical Poincar\'e group
it gives the standard relation $\ga^i\ga^j+\ga^j\ga^i=2g^{ij}\jb$.\\

\bd According to \cite{KG}, the condition 2 implies the
condition 1. Conversely, assume that the condition 1 holds.
Applying its both sides to $x^mx^n$, using \r{1.13}, \r{1.14}
and $Rg=g$, we get the condition 2. It remains to prove the
equivalence of conditions 2 and 3. 

Using \r{2.7a}, it is easy to check that the condition 2 is
equivalent  to
\be [\jb_4\ot(R+\jb)](W\ot\jb_4)W=2\cdot\jb_4\ot g. \e{3.2} \ee
In virtue of \r{2.8a}, \r{3.2} can be translated as 
\be [\jb_4\ot(R_{\lpi}+\jb)](N\ot\jb_4)N=2\cdot\jb_4\ot g_{\lpi}, \e{3.3}\ee
where
\[ R_{\lpi}=(\jb_2\ot X\ot\jb_2)(L\ot\te L)(\jb_2\ot
X^{-1}\ot\jb_2), \]
\[ g_{\lpi}=-2q^{1/2}(\jb_2\ot X\ot\jb_2)(E\ot\ta E). \]
This in turn means that
\be [\jb_2\ot(R_{\lpi}+\jb_{16})](N_1\ot\jb_4)N_2=2\cdot\jb_2\ot g_{\lpi}, 
\e{3.4}\ee
\be [\jb_2\ot(R_{\lpi}+\jb_{16})](N_2\ot\jb_4)N_1=2\cdot\jb_2\ot g_{\lpi}. 
\e{3.5}\ee
But we can write 
\[ N_1=2a(M\ot\jb_2),\qquad N_2=2bq^{-1/2}(\jb_2\ot
X^{-1}\ot T)(\ta E\ot E\ot\jb_2)\]
where $M:\cb\lra\cb^2\ot\cb^2$ and $T:\cb^2\ot\cb^2\lra\cb$ are
such that $M^{AB}=T_{AB}=\del^{A}{}_{B}$. Using this, 
\[ (\jb_2\ot T)(M\ot\jb_2)=(T\ot\jb_2)(\jb_2\ot M)=\jb_2 \]
and the 16 relations (2.1), (2.3)--(2.9), (2.17)--(2.20) and
(2.35)--(2.38) of \cite{POI}, after some computations we get
that the left hand sides of \r{3.4},\r{3.5} are both equal to
$2ab\jb_2\ot g_{\lpi}$. Therefore the condition 2 is equivalent to
$ab=1$.\ed

{\bf Remark.} Here we explain why each of conditions \r{3.4}--\r{3.5}
(expressing equalities of $32\ti 2$ matrices) leads to only one
numerical condition. Let us begin with \r{3.4}. In virtue of 
\r{2.11'},\r{2.12'} $(N_1\ot\jb_4)N_2$ 
intertwines $\wc$ with $z=\wc\ot\lpi\ot\lpi$ where $\lpi=w\ot\wb$.
Using Proposition 2.1 of \cite{POI}, 
\be \lpi\ot\lpi\ \simeq\ I\op w^1\op\ov{w^1}\op(w^1\ot\ov{w^1}). \e{3.5'}\ee
Thus 
\[ z\ \simeq\ w\op(w\op
w^{3/2})\op(w\ot\ov{w^1})\op(w\ot\ov{w^1}\op
w^{3/2}\ot\ov{w^1}). \]
Hence $(N_1\ot\jb_4)N_2$ can be nonzero only in first two
components of \r{3.5'}. Moreover, using the remarks after (2.25)
of \cite{POI}, $R_{\lpi}+\jb_{16}$ kills the second component of
\r{3.5'} and the LHS of \r{3.4} intertwines $\wc$ with
$\wc\ot I$ (here $I$ denotes the trivial subrepresentation of $\lpi\ot\lpi$). 
The same is true for the RHS (cf. \r{1.16}). Such intertwinners
are represented by numbers and \r{3.4} is equivalent to one
numerical condition. We treat \r{3.5} in a similar way. Our
remarks are valid in particular in the classical case.\\

The transformation
$a\mapsto x^{-1}a$, $b\mapsto xb$
($x\in\cb\backslash\{0\}$) corresponds to $\ga^i\mapsto
D\ga^i D^{-1}$ where $D=x\jb_2\op\jb_2$, which 
is equivalent to scaling of the undotted spinor.
Therefore we may set $a=b=1$ and obtain
\be \ga^i=\left( \ba{cc} 0&A_i\\ \si_i&0\ea\right)  \e{3.6}\ee
where $A_i$ are given by \r{2.12}.

In particular cases one gets:
\btr
\item[1)5)] $A_0=qt\si_0$, $A_1=-qt^{-1}\si_1$,
$A_2=-qt^{-1}\si_2$, $A_3=-qt\si_3$,
\item[2)6)] \[ A_0=q\left(\ba{cc}1-c^2&0\\0&1\ea\right),\quad 
A_1=-q\si_1,\quad
A_2=-q\si_2,\quad A_3=q\left(\ba{cc}-1-c^2&0\\0&1\ea\right),\]
\item[3)] \[ A_0=\left(\ba{cc}c^2(1-r)+1&c\\c&1\ea\right),\ 
A_1=\left(\ba{cc}-2c&-1\\-1&0\ea\right),\ 
A_2=\left(\ba{cc}0&i\\-i&0\ea\right),\ \]
\[ A_3=\left(\ba{cc}c^2(1-r)-1&c\\c&1\ea\right),\]
\item[4)]
\[ A_0=\left(\ba{cc}c^2+1&2c\\2c&1\ea\right),\quad
A_1=-\si_1,\quad A_2=-\si_2,\quad
A_3=\left(\ba{cc}c^2-1&2c\\2c&1\ea\right),\]
\item[7)] $A_0=(\sw-r)\si_0$, $A_1=(r+\sw)\si_1$, $A_2=(r-\sw)\si_2$,
$A_3=(r+\sw)\si_3$ 
\etr
($q=1$ for 1)--4) and $q=-1$ for 5)--7)).

In the following we shall study the analogues of certain
sesquilinear expressions which appear in the standard Dirac
theory. They involve the gamma matrix $\ga^0$. But in the
deformed case the corresponding matrix will be in general
different from $\ga^0$ and denoted by $A$. We set 
\be A=\left(\ba{cc}0&K^T\\K&0\ea\right), \e{3.16}\ee
where 
\be K=-E(E^{-1})^T\e{3.17}\ee
(as $2\ti 2$ matrices). For $\phi=\eps_a\ot\phi^a\in\te\cpi$ we
put $\phi^{\dag}=\lam_a\ot{\phi^a}^*$, $\bar\phi=\phi^{\dag}(A\ot I)$ where
$\lam_a\in(\cb^4)^*$ form a basis dual to the standard basis
$\eps_a$ of $\cb^4$. We shall prove that 
expressions like $\bar\phi\phi$ and
deformed Lagrangian $\lpi=\bar\phi(i\di-m)\phi$ transform
themselves in the same way as in the standard theory:

\bp \e{p3.2} For $\phi=\eps_a\ot\phi^a,
\chi=\eps_a\ot\chi^a\in\te\cpi$ one has
\be (\te\Psi\phi)^{\dag}(I\ot A\ot
I)(\te\Psi\chi)=\Psi[\phi^{\dag}(A\ot I)\chi], \e{3.17'}\ee
\be (\te\Psi\phi)^{\dag}(I\ot A\ot
I)[I\ot(i\di-m)](\te\Psi\chi)=\Psi[\phi^{\dag}(A\ot I)(i\di-m)\chi], 
\e{3.18}\ee
where $(b_a\ot c_a)^{\dag}=b_a^*\ot c_a^{\dag}$ for
$b_a\in\bpi$, $c_a\in\te\cpi$. \ep

\bd Using \r{2.2a} and \r{2.2'}, we notice that 
\[ wK(\wc)^T=wKE^Tw^T(E^{-1})^T=-wEw^T(E^{-1})^T=-E(E^{-1})^T=K.
\]
Applying the hermitian conjugation, one obtains
$\ov{\wc}K^T\wb^T=K^T$. These two equations and \r{2.2} give 
\be \bar \gpi A\gpi^T=A. \e{3.19}\ee
Therefore 
\[ (\te\Psi\phi)^{\dag}(I\ot A\ot
I)(\te\Psi\chi)=
[(\gpi_{a}{}^{l}{\phi^a}^{(1)})^*\ot\lam_l\ot{\phi^a}^{(2)*}](I\ot A\ot I)
[\gpi_{b}{}^{s}{\chi^b}^{(1)}\ot\eps_s\ot{\chi^b}^{(2)}] \]
\[ =
{\phi^a}^{(1)*}\gpi_{a}{}^{l}{}^*A_{ls}\gpi_{b}{}^{s}{\chi^b}^{(1)}
\ot{\phi^a}^{(2)*}
{\chi^b}^{(2)}
={\phi^a}^{(1)*}A_{ab}{\chi^b}^{(1)}\ot{\phi^a}^{(2)*}{\chi^b}^{(2)} \]
\[ = A_{ab}\Psi(\phi^a)^*\Psi(\chi^b)=\Psi(\lam_aA\eps_b\ot{\phi^a}^*\chi^b)=
\Psi(\phi^{\dag}(A\ot I)\chi) \]
and \r{3.17'} follows. Replacing $\chi$ by $(i\di-m)\chi$ and
using \r{2.4}, one gets \r{3.18}.\ed

{\bf Remark.} We were not able to get a similar fact for
$\bar\phi\ga^i\phi$. 

In particular cases $K$ equals:
\[ \mbox{ 1)2) }\ \  K=\jb_2,\quad \mbox{ 3)4) }\ \  K=\left(\ba{cc}1
&-2c\\0&1\ea\right),\quad \mbox{ 5)6)7) }\ \ K=-\jb_2. \]

\sectio{Solutions of the Dirac equation and momenta}

In this section we use \cite{KG} to get (in certain cases)
formal solutions of the Dirac equation. We introduce the momenta
for spin $1/2$ particles. They have good transformation
properties, commute with the Dirac operator, are
selfadjoint w.r.t. the inner product introduced in Section~2 and
in general differ from momenta for spin $0$ particles.
Unfortunately, only in some cases they are diagonalizable with
real eigenvalues. Considerations in this section are largely
formal (except of Proposition \re{p3.1}).

In the case when the Lorentz group is classical
(case 1), $t=1$), we get $R=\ta$. Then
we obtain formal solutions of the Dirac equation 
$(i\di-m)\phi=0$ ($m\geq0$) as in the case
2. of Section 4 of \cite{KG}. In these cases metric tensor 
$g=\diag(1,-1,-1,-1)$ and the gamma matrices
\[ \ga^0=\left(\ba{cc} 0&\jb_2\\ \jb_2&0\ea\right),\quad
\ga^i=\left(\ba{cc}0 &-\si_i\\\si_i&0\ea\right),\quad i=1,2,3, \]
are classical.
Solutions have form $\phi=v\ot e^{-ix^ap_a}$ (according to the
conventions of our paper, we use a different order of tensor
product than in \cite{KG}) where $p_a$ are real numbers and $v$
is a solution of
\be \ppi_j\ga^j v=mv \e{2.A}\ee
with
\[ m=\sqrt{\ppi_0^2-\ppi_1^2-\ppi_2^2-\ppi_3^2} \]
and $\ppi_j\in\rb$ obtained from $p_i$ as in \cite{KG}. We solve
\r{2.A} as in the standard theory:

for $m>0$ we get
\[ v=\left(\ba{c}\varphi\\m^{-1}\ppi_k\si_k\varphi\ea\right)\quad
(\varphi\in\cb^2), \]

for $m=0$ one has
\[ v\in\spa\{
\left(\ba{c}\ppi_1-i\ppi_2 \\-\ppi_0-\ppi_3\\0\\0\ea\right),\quad
\left(\ba{c}0\\0\\ \ppi_0+\ppi_3\\ \ppi_1+i\ppi_2\ea\right)\}\quad
(\ppi_0\neq-\ppi_3), \]
\[ v\in\spa\{
\left(\ba{c} 1\\0\\0\\0\ea\right),\quad \left(\ba{c} 0\\0\\0\\1\ea\right)\}
\quad (\ppi_0=-\ppi_3\neq0) \]

(or $v\in\cb^4$ if all $\ppi_i=0$ -- unphysical case).

Now let us consider the case 1. of Section~4 of \cite{KG} (for
$Z=0$). In addition to gamma matrices we need also
${}^*$-algebra $\fpi$ and its
${}^*$-representations $\pi$ in Hilbert spaces $H$ with bases
$e_k$, $k\in K$. We recall that $\fpi$ is generated by
$p^a$, $a=0,1,2,3$, satisfying $(p^a)^*=p^a$, $p^kp^l=R^{lk}{}_{ji}p^ip^j$.
Then solutions $\vhi$ of the Dirac equation $i\di\vhi=m\vhi$
($m\geq 0$) are provided in terms of
$\vhi_{sl}=(\id\ot\pi_{sl})(e^{-ix\ot p})$ where $x\ot p=x^a\ot g_{ab}p^b$,
$s,l\in K$, $\pi_{sl}(a)=(e_s\mid\pi(a)e_l)$ for $a\in\fpi$. 
Namely, $\vhi=\vhi_{vl}=\eps_i\ot v^{si}\vhi_{sl}$ ($i=1,2,3,4$)
for $v\in H\ot\cb^4$ such that $Uv=mv$. Here $U^{si}{}_{kn}=
\pi_{ks}(p^t)g_{at}(\ga^a)^i{}_n$, $m=\pi(g_{ij}p^jp^i)^{1/2}$
(we assume that it is a number, which holds e.g. for $\pi$ irreducible
as in \cite{KG}).

In the following we limit ourselves to the cases 1)-2) of
Section~1. Then $g_{ab}$ are real, $x\ot p$ is selfadjoint,
$e^{-ix\ot p}$ unitary and the components of 
$\vhi$ are bounded (since the summation over $s$ will be finite). 
Set $A=p^0+p^3$, $B=p^1-ip^2$, $B^*=p^1+ip^2$,
$D=p^0-p^3$. We have found the following admissible $\pi$:
\btr
\item[1)a)] $\pi_{abd}$ in $l^2(\zb)$ with orthonormal basis
$e_n$, $n\in\zb$, defined by
\[ \pi_{abd}(A)e_n=t^{-2n}ae_n, \]
\[ \pi_{abd}(B)e_n=be_{n-1}, \]
\[ \pi_{abd}(B^*)e_n=be_{n+1}, \]
\[ \pi_{abd}(D)e_n=t^{2n}de_n, \]
$a,d\in\rb$, $(a,d)\neq(0,0)$, $b>0$.
\item[b)] $\te\pi_{abd}$ in $\cb$ defined by
\[ \te\pi_{abd}(A)=a,\quad \te\pi_{abd}(B)=b,\quad
\te\pi_{abd}(B^*)=\bar b,\quad 
\te\pi_{abd}(D)=d, \]
$a,d\in\rb$, $b=0$ or $a=d=0$, $b\in\cbg$.
\item[2)a)] $\pi_{ad}$ in $l^2(\nb)$ with orthonormal basis
$e_n$, $n\in\nb$, defined by
\[ \pi_{ad}(A)e_n=c^2(a+nd)e_n, \]
\[ \pi_{ad}(B)e_n=cdn^{1/2}e_{n-1}, \]
\[ \pi_{ad}(B^*)e_n=cd(n+1)^{1/2}e_{n+1}, \]
\[ \pi_{ad}(D)e_n=de_n, \]
$d\in\rbg$, $a\in\rb$, $e_{-1}=0$.
\item[b)] $\te\pi_{ab}$ in $\cb$ defined by
\[ \te\pi_{ab}(A)=a,\quad \te\pi_{ab}(B)=b,\quad
\te\pi_{ab}(B^*)=\bar b,\quad 
\te\pi_{ab}(D)=0, \]
$a\in\rb$, $b\in\cb$.
\etr

After some computations one finds $U,v$ and finally solutions
$\vhi$ of the Dirac equation $i\di\vhi=m\vhi$
($m\geq0$). They are (up to
linear combinations) as follows:
\btr
\item[1)] $m=(t^{-1}ad-t\mid b\mid^2)^{1/2}$, 
\item[a)] $m>0$:
\[ \vhi=
\left(\ba{c}m\vhi_{nl} \\0 \\t^{2n-1}d\vhi_{nl} \\-tb\vhi_{n-1,l}
\ea \right),\quad
\left(\ba{c}0 \\m\vhi_{nl} \\-tb\vhi_{n+1,l}
\\t^{-2n-1}a\vhi_{nl}\ea \right); \]
$m=0$:
\[\vhi=
\left(\ba{c} b\vhi_{nl}\\dt^{2n-2}\vhi_{n-1,l} \\0 \\0 \ea\right), \quad
\left(\ba{c} 0 \\0 \\b\vhi_{nl} \\-at^{-2n}\vhi_{n-1,l} \ea\right).
\]
\item[b)] $m>0$ ($b=0$, $ad>0$):
\[ \vhi=
\left(\ba{c}m\vhi_{11} \\0 \\t^{-1}d\vhi_{11} \\0 \ea\right), \quad
\left(\ba{c}0 \\m\vhi_{11} \\0 \\t^{-1}a\vhi_{11} \ea\right); \]
$m=0$:
\[ \vhi= 
\left(\ba{c}0 \\ \vhi_{11}\\0 \\0 \ea\right),\ 
\left(\ba{c}0 \\0 \\ \vhi_{11}  \\0 \ea\right)\ \ (a=b=0),\]
\[ \vhi=\left(\ba{c}\vhi_{11} \\0 \\0 \\0 \ea\right), \
\left(\ba{c} 0\\0 \\0 \\ \vhi_{11} \ea\right)\ \  (d=b=0). \]
\item[2)a)] $m=c(ad)^{1/2}$, $m>0$:
\[ \vhi=
\left(\ba{c} m\vhi_{nl}\\0 \\d\vhi_{nl} \\-cdn^{1/2}\vhi_{n-1,l} 
\ea\right),\quad
\left(\ba{c} 0\\m\vhi_{nl}\\-cd(n+1)^{1/2}\vhi_{n+1,l}\\c^2(a+nd+d)\vhi_{nl}
\ea\right); \]
$m=0$:
\[ \vhi=
\left(\ba{c} c(n+1)^{1/2}\vhi_{n+1,l}\\ \vhi_{nl} \\0 \\0 \ea\right),\quad 
\left(\ba{c} 0 \\0 \\ \vhi_{nl} \\-cn^{1/2}\vhi_{n-1,l}
\ea\right). \]
\item[b)] $m=0$ for $b=0$, $a\neq0$:
\[ \vhi=
\left(\ba{c} \vhi_{11}\\0 \\0 \\0 \ea\right),\quad 
\left(\ba{c} 0\\0 \\0 \\ \vhi_{11} \ea\right) \] 
(we have omitted unphysical case $a=b=0$ when
$\vhi=\rho\ot\vhi_{11}$ with any $\rho\in\cb^4$).
\etr

Now we shall pass to the momenta for spin $1/2$ particles (in
general case).
Let us recall that for spin $0$ particles the momenta were
defined as $P_j=i\p_j$ and partial derivatives $\p_j$ can be
also obtained by $\p_j=(Y_j\ot\id)\Psi$
(see \r{1.15c}).
It suggests to define the momenta for spin $1/2$ particles as
\be \te P_j=i\te\p_j, \e{3.7}\ee
\be \te\p_j=(Y_j\ot\id)\te\Psi:\te\cpi\lra\te\cpi. \e{3.8}\ee
This choice is justified by the following

\bp\e{p3.1} Let us define $\te\p_j$ as in \r{3.8}. Then
\be \te\p_j\di=\di\te\p_j, \e{3.10}\ee
\be (\id\ot\te\p_j)\te\Psi(a)=(\la^{i}{}_{j}\ot
I)\te\Psi(\te\p_i(a)),\quad a\in\te\cpi.\e{3.11}\ee\ep 

{\bf Remark.} Proposition \re{p3.1} is valid for any homogeneous
quantum space endowed with the action of inhomogeneous quantum
group in the sense of \cite{KG} (with $\te F=0$). Equation \r{3.10} 
implies that
(up to technical difficulties) the deformed momenta are well
defined in the space of solutions of the Dirac equation
$(i\di-m)\phi=0$ ($m\geq0$). Equation \r{3.11} means that the
momenta for spin $1/2$ particles transform themselves in the
same way as the momenta for spin $0$ particles (cf. \r{1.15}).\\

\bd We set (cf. \r{1.15c})
\be \te\rho_{j}{}^{k}=
(X_{j}{}^{k}\ot\id)\te\Psi:\te\cpi\lra\te\cpi.\e{3.11''}\ee
Applying $Y_j\ot\id$ or $X_{j}{}^{k}\ot\id$ to \r{2.4}, we get
\r{3.10} and
\be \te\rho_{j}{}^{k}\di=\di\te\rho_{j}{}^{k}. \e{3.12}\ee

Next we are going to prove
\be Y_j*a=\la^{i}{}_{j}\{a*Y_i\}, \e{3.13}\ee
\be \{X_{k}{}^{j}*a\}\la^{i}{}_{j}=\la^{j}{}_{k}\{a*X_{j}{}^{i}\}, \e{3.14}\ee
$a\in\bpi$. For $a\in\api$ \r{3.13} is trivial while \r{3.14}
follows from \r{1.15a}, (1.5) of \cite{INH} and
invertibility of $S$. For $a=y^s$ \r{3.13} is trivial while
\r{3.14} follows from \r{1.15a},\r{1.15b} and (3.60) of
\cite{INH} ($p_i$ are called $y^i$ now). But 
the set of $a\in\bpi$ satisfying \r{3.13},\r{3.14} is an
algebra (use the homomorphism property of $\xb$ of \r{1.15b'}).
 Therefore \r{3.13},\r{3.14} follow for all $a\in\bpi$.

Moreover, using \r{2.1}, one obtains
\be (\id\ot\te\Psi)\te\Psi=(\de\ot\id)\te\Psi. \e{3.15}\ee
Now tensoring \r{3.13},\r{3.14} from right by $b\in\te\cpi$,
replacing $a\ot b$ by $\te\Psi(x)$, $x\in\te\cpi$, using
\r{3.15} and \r{3.8},\r{3.11''}, one obtains \r{3.11} and
\be (\id\ot\te\rho_{k}{}^{j})\te\Psi(a)(\la^{i}{}_{j}\ot I)=(\la^{j}{}_{k}\ot
I)\te\Psi(\te\rho_{j}{}^{i}(a)).\e{3.15'}\ee
In the same way, using $(\id\ot\Psi)\Psi=(\de\ot\id)\Psi$, one
can also get (1.11)--(1.12) of \cite{KG}.\ed

Using the properties of $\xb$ and \r{2.1}, one gets
\[ \te\p_m=X_m{}^j(\gpi_a{}^l)E_l{}^a\ot\p_j,\quad
\te\rho_m{}^k=X_m{}^j(\gpi_a{}^l)E_l{}^a\ot\rho_j{}^k \]
where $E_l{}^a$ are matrix units ($E_l{}^a\eps_b=\del^a{}_b\eps_l$).
First consider the case when the Lorentz group is classical [1),
$t=1$] and $s=1$. Then 
$f^{i}{}_{j}=\del^{i}{}_{j}\eps$. Using \r{1.15a}, one gets
\[ \te\p_m=\id\ot\p_m, \qquad  \te\rho_{m}{}^{k}=\id\ot\rho_{m}{}^{k}. \]
Taking the solution $\vhi=v\ot e^{-ix^ap_a}$ as in the beginning
of the present section, one obtains
\[ \te P_m\vhi=v\ot i\p_m e^{-ix^ap_a}=\ppi_m\vhi \]
(cf. the case 2. of Section~4 of \cite{KG}). Thus in this case
the momenta are just like the momenta for spin $0$ particles.
For all other cases it is easy to check that
$f^{i}{}_{m}(w^{A}{}_{B})\not\equiv\del^{i}{}_{m}\del^{A}{}_{B}$,
$X_{i}{}^{m}(\gpi_{a}{}^{l})\not\equiv\del^{m}{}_i\del^{l}{}_{a}$ and
\[ \te\p_i\not\equiv\id\ot\p_i, \qquad  
\te\rho_{i}{}^{k}\not\equiv\id\ot\rho_{i}{}^{k} \]
(acting on $\phi=\eps_a\ot x^b$ or $\phi=\eps_a\ot I$). It means
that in general momenta depend on spin.

Let us return to the general case. Set 
\[ F^t{}_r=g^{tm}X_m{}^j(\gpi_a{}^l)g_{jr}E_l{}^a, \]
$\te P^t=ig^{tm}\te\p_m$, $P^r=ig^{rj}\p_j$. Then
\[ \te P^t=F^t{}_r\ot P^r.\]
 The transformation property of the
momenta $\te P^t$ easily follows from \r{3.11}:
\[ (\id\ot\te P^t)\te\Psi(a)=[(g^{-1T}\la g^T)_m{}^t\ot
I]\te\Psi(\te P^m(a)). \]
Due to \r{3.10}, $\te P^t$ commute with $\di$. One can check that
$F^t{}_r$ are selfadjoint w.r.t. the inner product defined in
Section~2 ($\bar\phi F^t{}_r\psi=\overline{F^t{}_r\phi}\psi$ for
$\phi,\psi\in\cb^4$) while $P^r$ are selfadjoint according to
\cite{KG}. Therefore $\te P^t$ are selfadjoint. However, the
inner product of Section 2 is not positively defined and
selfadjointness doesn't guarantee diagonalizability or reality
of the spectrum as we shall see. Namely, after long computations
one gets the following form for $\te P^t$ in cases 1)-2), $Z=0$
(case 1. of Section 4 of \cite{KG} as considered above):

identifying $\vhi_{nl}$ with $e_n$ one can represent $\te P^t$
as operators in $\cb^4\ot H=(\cb^2\ot H)\oplus(\cb^2\ot H)$ such
that 
\[ \te P^t=\te R^t\oplus(\te R^t)^* \]
(this is related to the selfadjointness of $\te P^t$, $*$ is
given by the standard hermitian structures in $\cb^2$ and $H$)
where 
\btr
\item[1)]
\[ \te R^0=\frac s2\left(\ba{cc}t\pi(A)+t^{-1}\pi(D) & 0 \\
0  &t^{-1}\pi(A)+t\pi(D)  \ea\right),\]
\[ \te R^1=\frac s2\left(\ba{cc}t^{-1}\pi(B)^T+t\pi(B^*)^T &0  \\
0  &t^{-1}\pi(B^*)^T+t\pi(B)^T  \ea\right),\]
\[ \te R^2=i\frac s2\left(\ba{cc}t^{-1}\pi(B)^T -t\pi(B^*)^T & 0 \\
0  &t\pi(B)^T-t^{-1} \pi(B^*)^T  \ea\right),\]
\[ \te R^3=\frac s2\left(\ba{cc}t\pi(A)-t^{-1} \pi(D) &  0\\
0  &t^{-1}\pi(A)-t \pi(D)  \ea\right),\]
\item[2)]
\[ \te R^0=\frac s2\left(\ba{cc}\pi(A)+ \pi(D) &-c^2\pi(B)^T  \\
0  &\pi(A)+ \pi(D)  \ea\right),\]
\[ \te R^1=\frac s2\left(\ba{cc}\pi(B^*)^T+ \pi(B)^T & -c^2\pi(D) \\
0  &\pi(B^*)^T+ \pi(B)^T  \ea\right),\]
\[ \te R^2=i\frac s2\left(\ba{cc}\pi(B)^T- \pi(B^*)^T & c^2\pi(D) \\
0  &\pi(B)^T- \pi(B^*)^T  \ea\right),\]
\[ \te R^3=\frac s2\left(\ba{cc}\pi(A)- \pi(D) & -c^2\pi(B)^T \\
0  &\pi(A)- \pi(D)  \ea\right).\]
\etr

According to \r{3.10}, $\te P^t$ act in the subspace of solutions
of the Dirac equation. However, 
considering $\te P^1$, $\te P^2$ w.r.t. the found solutions,
there appears
complex spectrum (case 1)) or nondiagonalizability (case 2), $m>0$).
But in the case 2), $m=0$ all $\te P^r$ are diagonalizable with
real spectrum (in the subspace of found solutions). 
One can try to overcome this
difficulty by modifying the definition of $\te\p_m$. However,
one checks that it is not possible to 
get always (in all cases 1)2)) diagonalizability with real
spectrum 
(without spoiling the other properties of $\te P^m$) by
choosing another ansatz of the form $\te\p_m=M_m{}^j\ot\p_j$
($M_m{}^j\in\cb$, $M_m{}^j\not\equiv0$).

\begin{center} {\bf Acknowledgement} \end{center}

The  author is grateful to Professor W. Arveson and other faculty members
for their kind hospitality at UC Berkeley where a large part of this work
was done.
I am thankful to Professor J.~Wess, 
Dr M.~Klimek and Dr S.~Zakrzewski for interesting discussions.
\\\\

\end{document}